\documentclass[a4paper,11pt]{amsart}
\usepackage{graphicx}
\usepackage{mathptmx}
\usepackage{amsmath}
\usepackage{amssymb}
\usepackage{enumitem}
\usepackage{xcolor}
\usepackage{xparse}
\NewDocumentCommand{\eulerian}{omm}
 {%
  \genfrac<>{0pt}{}{#2}{#3}%
  \IfValueT{#1}{_{\!#1}}%
 }

 \newmuskip\pFqmuskip

\newcommand*\pFq[6][8]{%
  \begingroup 
  \pFqmuskip=#1mu\relax
  \mathchardef\normalcomma=\mathcode`,
  \mathcode`\,=\string"8000
  \begingroup\lccode`\~=`\,
  \lowercase{\endgroup\let~}\pFqcomma
  {}_{#2}F_{#3}{\left(\genfrac..{0pt}{}{#4}{#5}\bigg|#6\right)}%
  \endgroup
}
\newcommand{\pFqcomma}{{\normalcomma}\mskip\pFqmuskip}

\newtheorem{theorem}{Theorem}

\begin{document}

\title[Dimorphic Mersenne numbers and their applications]{Dimorphic Mersenne numbers and their applications}

\author{Taekyun  Kim}
\address{Department of Mathematics, Kwangwoon University, Seoul 139-701, Republic of Korea}
\email{tkkim@kw.ac.kr}

\author{DAE SAN KIM}
\address{Department of Mathematics, Sogang University, Seoul 121-742, Republic of Korea}
\email{dskim@sogang.ac.kr}


\begin{abstract}
The Mersenne primes are primes which can be written as some prime power of 2 minus 1. These primes were studied from antiquity in that their close connection with perfect numbers and even to present day in that their easiness for primality test.
In this paper, we introduce the dimorphic Mersenne numbers as a degenerate version of the Mersenne numbers and investigate some of their properties in connection with the degenerate Bernoulli polynomials and the incomplete Bell polynomials.
\end{abstract}

\maketitle

\section{Introduction}
In recent years, there have been active explorations for various degenerate versions of many special numbers and polynomials with diverse tools such as generating functions, combinatorial methods, $p$-adic analysis, umbral calculus, differential equations, probability theory, operator theory, analytic number theory and quantum physics. These were initiated by Carlitz in [3,4] and yielded many interesting results of arithmetical and combinatorial nature (see [7-9,12-14,16-18] and the references therein).\par
The aim of this paper is to introduce the dimorphic Mersenne numbers as a degenerate version of the Mersenne numbers and to find some of their applications in connection with the degenerate Bernoulli polynomials. The novelty of this paper is that it is the first paper which introduces the dimorphic Mersenne numbers and shows some of their applicatons associated with the degenerate Bernoulli numbers and the incomplete Bell polynomials. \par
The outline of this paper is as follows. In Section 1, we recall Mersenne numbers, the degenerate exponentials and the degenerate Bernoulli polynomials. We remind the reader of the degenerate Stirling numbers of the second kind, the incomplete Bell polynomials and the complete Bell polynomials. Section 2 is the main result of this paper. We first introduce the dimorphic Mersenne numbers as a degenerate version of the Mersenne numbers. Then, by using the generating function of the dimorphic Mersenne numbers in \eqref{16}, in Theorem 1, we derive a recurrence formula for the degenerate Bernoulli polynomials involving the dimorphic Mersenne numbers. We express a factor of the generating function of the dimorphic Mersenne numbers in terms of the incomplete Bell polynomials with arguments given by the degenerate Bernoulli polynomials (see \eqref{16}, \eqref{20}). In Theorem 2, we obtain from this expression a represnetation of the dimorphic Mersenne numbers in terms of those incomplete Bell polynomials with arguments given by the degenerate Bernoulli polynomials. Finally, in Theorem 3 we get an expression of the degenerate Bernoulli polynomials in terms of the same incomplete Bell polynomials with arguments given by the degenerate Bernoulli polynomials from Theorem 2. In the rest of this section, we recall the facts that are needed throughout this paper. \par
The Mersenne numbers are defined by 
\begin{displaymath}
	M_{n}=2^{n}-1,\quad (n\ge 0),\quad (\mathrm{see}\ [2]).
\end{displaymath}
$M_{p}$ is called a Mersenne prime if $p$ is prime and $M_{p}$ is also prime. 
For example, $M_{p}$ is a Mersenne prime for $p=2,3,5,7,13,17,19,31,61,89,107,127,521$ (see [2]). As of May 2022, the largest prime is the Mersenne prime $2^{82,589,933}-1$ which has 24,862,048 digits when it is written in base 10.
It is easy to show that 
\begin{equation}
\frac{z}{1-3z+2z^{2}}=\sum_{n=0}^{\infty}M_{n}z^{n},\quad (\mathrm{see}\ [2]). \label{1}
\end{equation} \par
For any nonzero $\lambda\in\mathbb{R}$, the degenerate exponentials are defined by 
\begin{equation}
e_{\lambda}^{x}(t)=(1+\lambda t)^{\frac{x}{\lambda}}=\sum_{n=0}^{\infty}(x)_{n,\lambda}\frac{t^{n}}{n!},\quad e_{\lambda}(t)=e_{\lambda}^{1}(t),\quad (\mathrm{see}\ [7-9,11-14]),\label{2}	
\end{equation}
where 
\begin{equation}
(x)_{0,\lambda}=1,\quad (x)_{n,\lambda}=x(x-\lambda)(x-2\lambda)\cdots(x-(n-1)\lambda),\quad (n\ge 1).\label{3}
\end{equation}
Note that 
\begin{displaymath}
	\lim_{\lambda\rightarrow 0}e_{\lambda}^{x}(t)=e^{xt},\quad \lim_{\lambda\rightarrow 0}(x)_{n,\lambda}=x^{n},\quad (n\ge 0). 
\end{displaymath} \par
In [3,4], Carlitz introduced the degenerate Bernoulli polynomials defined by  
\begin{equation}
\frac{t}{e_{\lambda}(t)-1}e_{\lambda}^{x}(t)=\sum_{n=0}^{\infty}\beta_{n,\lambda}(x)\frac{t^{n}}{n!}.\label{4}
\end{equation}
When $x=0$, $\beta_{n,\lambda}=\beta_{n,\lambda}(0),\ (n\ge 0)$, are called the degenerate Bernoulli numbers. Note that $\displaystyle \lim_{\lambda\rightarrow 0}\beta_{n,\lambda}=B_{n}\displaystyle$, where $B_{n}$ are the ordinary Bernoulli numbers, (see [1-19]). 
From \eqref{4}, we note that 
\begin{displaymath}
	\beta_{n,\lambda}(x)=\sum_{k=0}^{n}\binom{n}{k}(x)_{n-k,\lambda}\beta_{k,\lambda},\quad (n\ge 0),\quad (\mathrm{see}\ [3,4]). 
\end{displaymath} \par
It is known that the Stirling numbers of the second kind are defined by 
\begin{equation}
x^{n}=\sum_{k=0}^{n}S_{2}(n,k)(x)_{k},\quad (n\ge 0),\quad (\mathrm{see}\ [8]),\label{5}
\end{equation}
where $(x)_{0}=1,\ (x)_{n}=x(x-1)\cdots(x-n+1),\ (n\ge 1)$. 
From \eqref{5}, we note that 
\begin{equation}
\frac{1}{k!}\big(e^{t}-1\big)^{k}=\sum_{n=k}^{\infty}S_{2}(n,k)\frac{t^{n}}{n!},\quad (k \ge 0),\quad (\mathrm{see}\ [8]). \label{6}
\end{equation} \par
The incomplete Bell polynomials are defined by  
\begin{equation}
\frac{1}{k!}\bigg(\sum_{i=1}^{\infty}x_{i}\frac{t^{i}}{i!}\bigg)^{k}=\sum_{n=k}^{\infty}B_{n,k}(x_{1},\cdots,x_{n-k+1})\frac{t^{n}}{n!},\quad (k\ge 0),\quad (\mathrm{see}\ [1-8,10-13,15]),\label{7}
\end{equation}
More explicitly, they are given by
\begin{equation}
\begin{aligned}
&B_{n,k}(x_{1},x_{2},\dots,x_{n-k+1})\\
&=\sum_{\substack{l_{1}+\cdots+l_{n-k+1}=k\\ l_{1}+2l_{2}+\cdots+(n-k+1)l_{n-k+1}=n}}\frac{n!}{l_{1}!l_{2}!\cdots l_{n-k+1}!}	\bigg(\frac{x_{1}}{1!}\bigg)^{l_{1}} \bigg(\frac{x_{2}}{2!}\bigg)^{l_{2}}\cdots \bigg(\frac{x_{n-k+1}}{(n-k+1)!}\bigg)^{l_{n-k+1}},
\end{aligned}\label{8}
\end{equation}
where the sum runs over all nonnegetive integers $l_{1}, \dots, l_{n-k+1}$ satisfying $l_{1}+\cdots+l_{n-k+1}=k$ and $l_{1}+2l_{2}+\cdots+(n-k+1)l_{n-k+1}=n$. \par
The complete Bell polynomials are given by 
\begin{align}
\exp\bigg(\sum_{i=1}^{\infty}x_{i}\frac{t^{i}}{i!}\bigg)&=\sum_{n=0}^{\infty}B_{n}(x_{1},\dots,x_{n})\frac{t^{n}}{n!} \label{9}\\
&=1+\sum_{n=1}^{\infty}B_{n}(x_{1},x_{2},\dots,x_{n})\frac{t^{n}}{n!},\quad (\mathrm{see}\ [1-8,10-13,15]). \nonumber	
\end{align}
Then we have
\begin{displaymath}
B_{0}(x_{1},\dots,x_{n})=1,\ B_{n}(x_{1},\dots,x_{n})=\sum_{l_{1}+2l_{2}+\cdots+nl_{n}=n}\frac{n!}{l_{1}!\cdots l_{n}!}\bigg(\frac{x_{1}}{1!}\bigg)^{l_{1}}\cdots\bigg(\frac{x_{n}}{n!}\bigg)^{l_{n}},\quad (n \ge 1),
\end{displaymath}
where the sum runs over all nonnegetive integers $l_{1}, \dots, l_{n}$ satisfying $l_{1}+2l_{2}+\cdots+nl_{n}=n$. \par
From \eqref{7} and \eqref{8}, we note that 
\begin{align}
\exp\bigg(\sum_{i=1}^{\infty}x_{i}\frac{t^{i}}{i!}\bigg) 
&=1+\sum_{k=1}^{\infty}\frac{1}{k!}\bigg(\sum_{i=1}^{\infty}x_{i}\frac{t^{i}}{i!}\bigg)^{k} \label{10}\\ 
&=1+\sum_{k=1}^{\infty}\sum_{n=k}^{\infty}B_{n,k}(x_{1},x_{2},\dots,x_{n-k+1})\frac{t^{n}}{n!} \nonumber \\
&=1+\sum_{n=1}^{\infty}\sum_{k=1}^{n}B_{n,k}(x_{1},\dots,x_{n-k+1})\frac{t^{n}}{n!}.\nonumber 	
\end{align}
By \eqref{9} and \eqref{10}, we get
\begin{equation}
B_{n}(x_{1},x_{2},\dots,x_{n})=\sum_{k=1}^{n}B_{n,k}(x_{1},\dots,x_{n-k+1}),\quad (n\ge 1). \label{11}
\end{equation}
It is known that the Bell polynomials are given by 
\begin{displaymath}
\phi_{n}(x)=\sum_{k=0}^{n}S_{2}(n,k)x^{k},\quad (n\ge 0),\quad (\mathrm{see}\ [16,17]).
\end{displaymath}
From \eqref{6}, \eqref{7} and \eqref{11}, we note that 
\begin{displaymath}
B_{n,k}(1,1,\dots,1)=S_{2}(n,k),\quad B_{n}(x,x,\dots,x)=\phi_{n}(x),\quad (n\ge 0). 
\end{displaymath}

\section{Dimorphic Mersenne numbers and their applications}
Now, we consider the dimorphic Mersenne numbers given by 
\begin{equation}
M_{n,\lambda}=(2)_{n,\lambda}-(1)_{n,\lambda},\quad (\lambda\in\mathbb{R},\quad n\ge 0). \label{12}	
\end{equation}
Note that $\lim_{\lambda\rightarrow 0}M_{n,\lambda}=M_{n}$. 
We observe that
\begin{align}
	e_{\lambda}^{2}(t)-e_{\lambda}(t)&=\sum_{n=0}^{\infty}(2)_{n,\lambda}\frac{t^{n}}{n!}-\sum_{n=0}^{\infty}(1)_{n,\lambda}\frac{t^{n}}{n!} \label{13} \\
	&=\sum_{n=0}^{\infty}\Big((2)_{n,\lambda}-(1)_{n,\lambda}\Big)\frac{t^{n}}{n!}. \nonumber
\end{align}
Thus, by \eqref{12} and \eqref{13}, we have 
\begin{equation}
e_{\lambda}^{2}(t)-e_{\lambda}(t)=\sum_{n=0}^{\infty}M_{n,\lambda}\frac{t^{n}}{n!}. \label{14}
\end{equation}
Thus, by \eqref{14} and noting that $M_{0,\lambda}=0$, we see that 
\begin{equation}
\frac{1}{t}\big(e_{\lambda}^{2}(t)-e_{\lambda}(t)\big)=\frac{1}{t}\sum_{n=1}^{\infty}M_{n,\lambda}\frac{t^{n}}{n!}=\sum_{n=0}^{\infty}\frac{M_{n+1,\lambda}}{n+1}\frac{t^{n}}{n!}. \label{15}
\end{equation}
From \eqref{15}, we have 
\begin{align}
\sum_{n=0}^{\infty}\frac{M_{n+1,\lambda}}{n+1}\frac{t^{n}}{n!}&=\frac{1}{t}\big(e_{\lambda}^{2}(t)-e_{\lambda}(t)\big)	=\frac{e_{\lambda}(t)}{t}\big(e_{\lambda}(t)-1\big)\label{16} \\
&=\frac{e_{\lambda}^{x+1}(t)}{te_{\lambda}^{x}(t)}(e_{\lambda}(t)-1). \nonumber
\end{align} \par
Now, we observe that 
\begin{align}
e_{\lambda}^{x+1}(t)&=\frac{t}{e_{\lambda}(t)-1}e_{\lambda}^{x}(t) \frac{e_{\lambda}^{x+1}(t)}{te_{\lambda}^{x}(t)}\big(e_{\lambda}(t)-1\big) \label{17}\\
&=\sum_{l=0}^{\infty}\beta_{l,\lambda}(x)\frac{t^{l}}{l!}\sum_{m=0}^{\infty}\frac{M_{m+1,\lambda}}{m+1}\frac{t^{m}}{m!}\nonumber\\
&=\sum_{n=0}^{\infty}\sum_{l=0}^{n}\binom{n}{l}\beta_{l,\lambda}(x)\frac{M_{n-l+1,\lambda}}{n-l+1}\frac{t^{n}}{n!}. \nonumber
\end{align}
On the other hand, by \eqref{2}, we also have
\begin{equation}
e_{\lambda}^{x+1}(t)=\sum_{n=0}^{\infty}(x+1)_{n,\lambda}\frac{t^{n}}{n!}.\label{18}	
\end{equation}
From \eqref{17} and \eqref{18}, we note that 
\begin{align}
(x+1)_{n,\lambda} &=\sum_{l=0}^{n}\binom{n}{l}\beta_{l,\lambda}(x)\frac{M_{n-l+1,\lambda}}{n-l+1}\label{19}\\
&=\beta_{n,\lambda}(x)+\sum_{l=0}^{n-1}\binom{n}{l}\beta_{l,\lambda}(x)\frac{M_{n-l+1,\lambda}}{n-l+1}.\nonumber	
\end{align}
Therefore, by \eqref{19}, we obtain the following theorem.
\begin{theorem}
For $n\ge 0$, we have 
\begin{displaymath}
\beta_{n,\lambda}(x)=(x+1)_{n,\lambda}-\sum_{l=0}^{n-1}\binom{n}{l}\beta_{l,\lambda}(x)\frac{M_{n-l+1,\lambda}}{n-l+1}. 
\end{displaymath}
\end{theorem}
From \eqref{7}, we note that 
\begin{align}
\frac{e_{\lambda}(t)-1}{te_{\lambda}^{x}(t)}&=\frac{1}{\frac{t}{e_{\lambda}(t)-1} e_{\lambda}^{x}(t)}=\frac{1}{1+\frac{t}{e_{\lambda}(t)-1}e_{\lambda}^{x}(t)-1}\label{20} \\
&=\frac{1}{1+\sum_{n=1}^{\infty}\beta_{n,\lambda}(x)\frac{t^{n}}{n!}}=\sum_{k=0}^{\infty}(-1)^{k}\bigg(\sum_{i=1}^{\infty}\beta_{i,\lambda}(x)\frac{t^{i}}{i!}\bigg)^{k} \nonumber \\
&=1+\sum_{k=1}^{\infty}(-1)^{k}k!\frac{1}{k!}\bigg(\sum_{i=1}^{\infty}\beta_{i,\lambda}(x)\frac{t^{i}}{i!}\bigg)^{k}\nonumber \\
&=\sum_{k=1}^{\infty}(-1)^{k}k!\sum_{n=k}^{\infty}B_{n,k}\big(\beta_{1,\lambda}(x),\beta_{2,\lambda}(x),\dots,\beta_{n-k+1,\lambda}(x)\big)\frac{t^{n}}{n!}+1 \nonumber \\
&=\sum_{n=1}^{\infty}\sum_{k=1}^{n}(-1)^{k}k!B_{n,k}\big(\beta_{1,\lambda}(x),\beta_{2,\lambda}(x),\dots,\beta_{n-k+1,\lambda}(x)\big)\frac{t^{n}}{n!}+1.\nonumber
\end{align}\\
On the other hand, by \eqref{16}, we get 
\begin{align}
&1+\sum_{n=1}^{\infty}\frac{M_{n+1,\lambda}}{n+1}\frac{t^{n}}{n!}=\sum_{n=0}^{\infty}\frac{M_{n+1,\lambda}}{n+1}\frac{t^{n}}{n!}=\frac{e_{\lambda}^{x+1}(t)}{te_{\lambda}^{x}(t)}(e_{\lambda}(t)-1)\label{21}\\
&=e_{\lambda}^{x+1}(t)\bigg(1+\sum_{j=1}^{\infty}\sum_{k=1}^{j}(-1)^{k}k!B_{j,k}\big(\beta_{1,\lambda}(x),\dots,\beta_{j-k+1,\lambda}(x)\big)\frac{t^{j}}{j!}\bigg)\nonumber \\
&=\sum_{l=0}^{\infty}(x+1)_{l,\lambda}\frac{t^{l}}{l!}\bigg(1+\sum_{j=1}^{\infty}\sum_{k=1}^{j}(-1)^{k}k!B_{j,k}\big(\beta_{1,\lambda}(x),\dots,\beta_{j-k+1,\lambda}(x)\big)\frac{t^{j}}{j!}\bigg)\nonumber \\
&=\sum_{n=0}^{\infty}(x+1)_{n,\lambda}\frac{t^{n}}{n!}+\sum_{n=1}^{\infty}\sum_{j=1}^{n}\sum_{k=1}^{j}\binom{n}{j}(x+1)_{n-j,\lambda}(-1)^{k}k!B_{j,k}\big(\beta_{1,\lambda}(x),\dots,\beta_{j-k+1,\lambda}(x)\big)\frac{t^{n}}{n!} \nonumber \\
&=1+\sum_{n=1}^{\infty}\bigg((x+1)_{n,\lambda}+\sum_{j=1}^{n}\sum_{k=1}^{j}\binom{n}{j}(x+1)_{n-j,\lambda}(-1)^{k}k!B_{j,k}\big(\beta_{1,\lambda}(x),\dots,\beta_{j-k+1,\lambda}(x)\big)\bigg)\frac{t^{n}}{n!}. \nonumber
\end{align}
Therefore, by comparing the coefficients on both sides of \eqref{21}, we obtain the following theorem. 
\begin{theorem}
For $n\in\mathbb{N}$ and any $x$, we have 
\begin{displaymath}
\frac{M_{n+1,\lambda}}{n+1}=(x+1)_{n,\lambda}+\sum_{i=1}^{n}\sum_{k=1}^{j}\binom{n}{j}(x+1)_{n-j,\lambda}(-1)^{k}k!B_{j,k}\big(\beta_{1,\lambda}(x),\dots,\beta_{j-k+1,\lambda}(x)\big). 
\end{displaymath}
\end{theorem}
We note that 
\begin{align}
B_{n,1}\big(\beta_{1,\lambda}(x),\beta_{2,\lambda}(x),\dots,\beta_{n,\lambda}(x)\big)&=\sum_{\substack{l_{1}+\cdots+l_{n}=1\\ l_{1}+2l_{2}+\cdots+nl_{n}=n}}\frac{n!}{l_{1}!l_{2}!\cdots l_{n}!}\bigg(\frac{\beta_{1,\lambda}(x)}{1!}\bigg)^{l_{1}}\cdots \bigg(\frac{\beta_{n,\lambda}(x)}{n!}\bigg)^{l_{n}}\label{22}\\
&=n!\frac{\beta_{n,\lambda}(x)}{n!}=\beta_{n,\lambda}(x). \nonumber	
\end{align}
From Theorem 2 and \eqref{22}, we further observe that
\begin{align}
&\frac{M_{n+1,\lambda}}{n+1}=(x+1)_{n,\lambda}+\sum_{j=1}^{n}\sum_{k=1}^{j}\binom{n}{j}(x+1)_{n-j,\lambda}(-1)^{k}k!B_{j,k}\big(\beta_{1,\lambda}(x),\beta_{2,\lambda}(x),\dots,\beta_{j-k+1,\lambda}(x)\big)\label{23} \\
&=(x+1)_{n,\lambda}+\sum_{k=1}^{n}(-1)^{k}k!B_{n,k}\big(\beta_{1,\lambda}(x),\dots,\beta_{n-k+1,\lambda}(x)\big)\nonumber \\
&\quad +\sum_{j=1}^{n-1}\sum_{k=1}^{j}\binom{n}{j}(x+1)_{n-j,\lambda}(-1)^{k}k!B_{j,k}\big(\beta_{1,\lambda}(x),\dots,\beta_{j-k+1,\lambda}(x)\big) \nonumber \\
&=(x+1)_{n,\lambda}-B_{n,1}\big(\beta_{1,\lambda}(x),\dots,\beta_{n,\lambda}(x)\big)+\sum_{k=2}^{n}(-1)^{k}k!B_{n,k}\big(\beta_{1,\lambda}(x), \beta_{2,\lambda}(x),\dots, \beta_{n-k+1,\lambda}(x)\big)\nonumber \\
&\quad +\sum_{j=1}^{n-1}\sum_{k=1}^{j}\binom{n}{j}(x+1)_{n-j,\lambda}(-1)^{k}k!B_{j,k}\big(\beta_{1,\lambda}(x), \beta_{2,\lambda}(x),\dots, \beta_{j-k+1,\lambda}(x)\big) \nonumber \\
&= -\beta_{n,\lambda}(x)+(x+1)_{n,\lambda}+\sum_{k=2}^{n}(-1)^{k}k!B_{n,k}\big(\beta_{1,\lambda}(x), \beta_{2,\lambda}(x),\dots, \beta_{n-k+1,\lambda}(x)\big) \nonumber \\
&\quad +\sum_{j=1}^{n-1}\sum_{k=1}^{j}\binom{n}{j}(x+1)_{n-j,\lambda}(-1)^{k}k!B_{j,k}\big(\beta_{1,\lambda}(x), \beta_{2,\lambda}(x),\dots, \beta_{j-k+1,\lambda}(x)\big).\nonumber
\end{align}
Thus, by \eqref{23}, we get 
\begin{equation}
\begin{aligned}
\beta_{n,\lambda}(x)&=(x+1)_{n,\lambda}-\frac{M_{n+1,\lambda}}{n+1}+\sum_{k=2}^{n}(-1)^{k}k!B_{n,k}\big(\beta_{1,\lambda}(x),\beta_{2,\lambda}(x),\dots,\beta_{n-k+1,\lambda}(x)\big)	\\
&\quad +\sum_{j=1}^{n-1}\sum_{k=1}^{j}\binom{n}{j}(x+1)_{n-j,\lambda}(-1)^{k}k!B_{j,k}\big(\beta_{1,\lambda}(x),\beta_{2,\lambda}(x),\dots,\beta_{j-k+1,\lambda}(x)\big).
\end{aligned}\label{24}
\end{equation}
Therefore, by \eqref{24}, we obtain the following theorem. 
\begin{theorem}
For $n\in\mathbb{N}$, we have 
\begin{align*}
\beta_{n,\lambda}(x)&=(x+1)_{n,\lambda}-\frac{M_{n+1}}{n+1}+\sum_{k=2}^{n}(-1)^{k}k!B_{n,k}\big(\beta_{1,\lambda}(x),\beta_{2,\lambda}(x),\dots,\beta_{n-k+1,\lambda}(x)\big)	\\
&\quad +\sum_{j=1}^{n-1}\sum_{k=1}^{j}\binom{n}{j}(x+1)_{n-j,\lambda}(-1)^{k}k!B_{j,k}\big(\beta_{1,\lambda}(x),\beta_{2,\lambda}(x),\dots,\beta_{j-k+1,\lambda}(x)\big).
\end{align*}
\end{theorem}

\section{Conclusion}
Study of various degenerate versions of some special numbers and polynomials, which began with the papers [3,4] by Carlitz, regained recent interests of some mathematicians and led to unexpected introduction of degenerate gamma functions and degenerate umbral calculus (see [9,14]).
\indent In this paper, the dimorphic Mersenne numbers were introduced as a degenerate version of the Mersenne numbers and some of their properties was investigated in connection with the degenerate Bernoulli polynomials and the incomplete Bell polynomials. \par
It is one of our future research projects to continue to study various degenerate versions of some special numbers and polynomials and to find their applications to physics, science and engineering as well as to mathematics.

\end{document}